\documentclass[11pt,reqno,oneside]{amsart}

\usepackage[utf8]{inputenc}
\usepackage[T1]{fontenc}

\usepackage{amsmath,amsthm,amssymb}

\usepackage{graphics}
\usepackage[pagebackref]{hyperref}
\renewcommand*\backref[1]{\ifx#1\relax \else (Cited on page #1) \fi}

\usepackage{amsmath,amsthm,amssymb}

\usepackage{graphics}
\usepackage{hyperref}
\usepackage[dvipsnames]{xcolor}

\usepackage{booktabs}

\usepackage{mathtools}
\mathtoolsset{showonlyrefs}

\usepackage{enumitem}
\newlist{todolist}{itemize}{2}
\setlist[todolist]{label=$\square$}
\usepackage{pifont}
\usepackage[square,sort,comma,numbers]{natbib}

\definecolor{darkblue}{rgb}{0.0,0.0,0.3}
\hypersetup{colorlinks,breaklinks,
  linkcolor=darkblue,urlcolor=darkblue,
anchorcolor=darkblue,citecolor=darkblue}

\theoremstyle{plain}
\newtheorem*{theorem*}{Theorem}

\newtheorem*{proposition*}{Proposition}

\newtheorem*{corollary*}{Corollary}

\theoremstyle{definition}

\newtheorem*{algorithm*}{Algorithm}

\newtheorem*{experiment*}{Experiment}

\DeclareMathOperator{\SL}{SL}

\title{A database of rigorous Maass forms}
\author{David Lowry-Duda}
\date{}

\begin{document}

\begin{abstract}
  We announce a database of rigorously computed Maass forms on congruence
  subgroups $\Gamma_0(N)$ and briefly describe the methods of computation.
\end{abstract}

\maketitle

\section{Computation of Maass Forms}

Let $\Gamma_0(N)$ be a congruence subgroup in $\SL(2, \mathbb{Z})$.
A matrix $\gamma = \left( \begin{smallmatrix} a & b \\ c & d \end{smallmatrix}
\right)$ acts on the upper halfplane $\mathcal{H}$ in the standard way, given
by $\gamma z = \frac{az + b}{cz + d}$.
The quotient space $\Gamma_0(N) \backslash \mathcal{H}$ is a noncompact Riemann
surface.
In the hyperbolic metric on $\mathcal{H}$, the Laplace-Beltrami operator,
$\Delta$, takes the form
\begin{equation}
  \Delta
  =
  y^2 \Big(
    \frac{\partial^2}{\partial x^2} + \frac{\partial^2}{\partial y^2}
  \Big).
\end{equation}
Maass cuspforms of weight $0$ on $\Gamma_0(N)$ are eigenfunctions $f$ of the
Laplace-Beltrami operator satisfying
\begin{equation}
  \Delta f + \lambda f = 0, \quad
  f(\gamma z) = f(z) \text{ for all } \gamma \in \Gamma_0(N),
  \quad \text{and} \quad
  f \in L^2(\Gamma_0(N) \backslash \mathcal{H}).
\end{equation}
Maass cuspforms form the discrete component of the spectral resolution of $\Delta$
and therefore form building blocks for all $L^2(\Gamma_0(N) \backslash
\mathcal{H})$.

Despite this, it is difficult to give any nontrivial example of a Maass
cuspform.
The problem of giving numerical examples of Maass cusp forms has been widely
considered since the 1970s.
In~\cite{then2005maass}, H.\ Then gives an extensive list of references to
earlier computations of Maass forms.
Until 2006, methods to generate numerical examples were
fundamentally \emph{heuristic}.

Booker, Str\"ombergsson and Venkatesh~\cite{booker2006effective} gave the first
method of computing \emph{rigorous} Maass forms on $\SL(2, \mathbb{Z})$.
Stated more precisely, they gave a method prove bounds for how close a putative
Maass form is to an actual Maass form.
There are now three algorithms to generate rigorous Maass cuspforms:
\begin{enumerate}
  \item \emph{Quasimode construction}, a generalization
  of~\cite{booker2006effective} to general level and character by
  Child~\cite{child22}.
  \item A rigorous implementation of the \emph{Selberg trace formula} due to
  Seymour-Howell~\cite{seymour2022rigorous},
  \item A rigorous version of \emph{Hejhal's algorithm} due to Seymour-Howell
  and Lowry-Duda, to be described in forthcoming work~\cite{ldshHejhal}.
\end{enumerate}
It is notable that Booker (and to a lesser extent, Str\"ombergsson) helped
develop each of these algorithms.

In this note, we describe how an initial database of 35416 rigorous Maass
cuspforms (of weight $0$, on congruence subgroups $\Gamma_0(N)$ with $N$
squarefree) was computed using these algorithms and inserted into the
LMFDB~\cite{lmfdb}.
We briefly describe the three algorithms and how they interact in \S2.
In \S3, we detail the data now available.
And in \S4, we give additional comments on the current and future status of
this database.

\section{Algorithms for Computing Maass Forms}\label{sec:algorithms}

Each Maass form considered here can be written in the form
\begin{equation}\label{eq:maass_expansion}
  f_j(z)
  =
  \sum_{n \neq 0}
  a_j(n) \sqrt{y} K_{i{r_j}}(2 \pi \lvert n \rvert y) e^{2 \pi i n x}
\end{equation}
where $z = x + iy$, $K_\nu(y)$ is the modified Bessel function of the second
kind, and the eigenvalue $\lambda$ of $f_j$ is given by $\lambda = \frac{1}{4}
+ r_j^2$ for $r_j \in \mathbb{R}$.
``Rigorously computing'' a Maass cusp form here means giving rigorous
approximations to the eigenvalue $\lambda$ and the first several coefficients
$a_j(n)$.
Each of the three methods for computation have different strengths.

\subsection{Rigorous Trace Formula}

With extreme simplification, the Selberg trace formula in~\cite{seymour2022rigorous}
is a formula for
\begin{equation}\label{eq:tracelhs}
  \sum_{j > 0} F(r_j) a_j(n),
\end{equation}
where the sum is over all Maass forms $f_j$ on $\Gamma_0(N)$.
The formula gives~\eqref{eq:tracelhs} in terms of a complicated sum of
conjugacy classes of the group $\Gamma_0(N)$ weighted by a function $G$, where
$F$ and $G$ are test functions related by
\begin{equation}
  G(u) = \int_{-\infty}^\infty F(r) e^{- 2 \pi i r u} dr.
\end{equation}
By carefully making different choices of test functions, it's possible to
isolate individual coefficients and eigenvalues.
One major difficulty in this approach is that class numbers of quadratic fields
naturally appear but are hard to compute.
This limits the range of computation (though this relationship can
be used to compute \emph{more} class numbers~\cite{booker2024unconditional}).

The most important benefit of this algorithm is that it can \emph{guarantee
that all Maass forms have been found} in an eigenvalue range.
In practice, one can try to use the trace formula to give initial
approximations to Maass forms and use other methods to refine.
On the other hand, the current implementation is currently restricted to
squarefree level and weight $0$.

\subsection{Quasimode construction}

In~\cite{child22}, Child extends the quasimode construction technique
in~\cite{booker2006effective} to general level.
If $\widetilde{f}_j$ is a putative Maass form and $(\Delta - \lambda)
\widetilde{f}_j$ has small $L^2$ norm, then a spectral resolution shows that
$\widetilde{f}_j$ is close to a true eigenfunction.
Child (and BSV) show that it is sufficient to obtain strong bounds along the
boundary of the fundamental domain.
This method can be used to certify precise (but heuristic) approximations to
Maass forms.

\subsection{Rigorous Hejhal}

Hejhal's algorithm~\cite{hejhal1999eigenfunctions} is in practice one of the
better algorithms for producing heuristic approximations to Maass forms.
Str\"omberg described how to adapt Hejhal's heuristic algorithm for general
level (and even general multiplier system)
in~\cite{stromberg2005computational}.

The fundamental idea is to use automorphy to construct an approximation linear
system of equations for the coefficients.
Rapid decay from the Bessel functions in~\eqref{eq:maass_expansion} shows that
the truncation $\widetilde{f}_j$ to the first $M$ coefficients is close to
$f_j$.
Taking an appropriate linear combination of $\widetilde{f}_j$ at points $z_j =
x_j + iY$ along a fixed horocycle, one obtains equations of the form
\begin{equation}
  a_j(n) \sqrt{Y} K_{ir_j}(2 \pi \lvert n \rvert Y)
  =
  \frac{1}{2Q} \sum_{j = 1 - Q}^Q \widetilde{f}(z_j) e(-nx_j)
  +
  (\text{truncation error}).
\end{equation}
If the horocycle points are chosen so that each $z_j$ is outside the
fundamental domain, then we compute the pullbacks $z_j^*$ to the fundamental
domain and insert those instead.
This introduces nonlinear mixing into the system, and solving gives approximate
Maass forms.
Of course, we don't actually know $r_j$.
Instead we guess a value and obtain heuristic coefficients.

The rigorous implementation of Hejhal's algorithm makes the truncation and
other errors explicit and tracks how an initial approximation (e.g.\ coming
from the trace formula), guaranteed to some initial precision, behaves under
iteration of Hejhal's algorithm.
When the initial approximation is sufficiently accurate, rigorously applying
Hejhal's algorithm refines and produces provably better approximations.
Unfortunately, this requires initial approximations; and if the initial
approximations are not sufficiently strong, this algorithm may fail to improve
the precision.

\section{Data Computed}\label{sec:data}

We computed $35416$ rigorously Maass cusp forms on $\Gamma_0(N)$ for each
squarefree $N$ from $1$ to $105$.
For each Maass form, we compute the eigenvalue, the first $1000$ coefficients,
and a portrait (constructed using similar methods as~\cite{lowryduda_visualizing}).
As the eigenvalue and coefficients are approximations to (a priori,
transcendental) real numbers, they are stored as pairs of $(\text{center},
\text{error})$.
In practice, \texttt{arb}~\cite{arb} was used for both data generation and
processing.

The database consists of approximately $4.954$ GB of Maass form data in total.

Level $1$ is distinguished, because both trace formula methods and verification
methods work much better there.
There are $2202$ Maass forms of level $1$ in the database, each verified using
the quasimode construction above.

For each level $2 \leq n \leq 105$, we computed as many Maass forms as we could
while guaranteeing that no eigenvalues were omitted.
The limiting factor is the trace formula, which is our only way to be certain
that every Maass form was detected.

\section{Comments on Database Construction}\label{sec:comments}

Finally, we conclude with several small comments.

\begin{enumerate}
  \item
  It is often easier to produce good estimates to the eigenvalues than to
  the coefficients.
  The current database includes $15423$ forms with coefficients that are too
  imprecise to detect the Fricke sign.
  Different methods to handle these problems are in development.

  \item
  By restricting to squarefree level, we omit Maass forms coming from induced
  representations of Hecke characters.
  These are Maass cuspforms that are \emph{explicitly computable} (see
  e.g.~\cite{maass1949automorphe}).
  Algorithms to compute these efficiently have been recently implemented in
  PARI/GP by Molin and Page~\cite{molin2022computing}.

  We also miss more complicated (and potentially interesting) Artin
  representations.

  Forthcoming work of Booker, Bober, Knightly, Krishnamurthy, Lee, Lowry-Duda,
  and Seymour-Howell seeks to work out an explicit, computable trace formula
  for general level and weight.

  \item
  The database currently omits all $L$-functions of Maass forms.
  In principle, these can be constructed from current data.
  But in practice existing algorithms and implementations are designed to work
  with algebraic $L$-functions.
  On the other hand, we shouldn't expect $L$-functions of an individual Maass
  form to appear elsewhere in the LMFDB.\@
\end{enumerate}

\section*{Acknowledgements}

We thank Andrew Booker, Min Lee, Brendan Hassett, Andrei Seymour-Howell, and
Drew Sutherland for guidance and support throughout the computation.
We also thank John Voight for several suggestions on how to improve the initial
form of the database, as well as David Roe, who helped improve the database
following these suggestions.
This task was easier thanks to the work of Edgar Costa and David Roe, and the
entire process they've helped create around contributing to the
LMFDB~\cite{costa2021zen}.

This work was supported by the Simons Collaboration in Arithmetic
Geometry, Number Theory, and Computation via the Simons Foundation grant
546235.

\bibliographystyle{alpha}
\bibliography{bibfile}

\end{document}